\documentclass{gtart}

\def\ifplaintex{\expandafter\ifx\csname documentclass\endcsname\relax}

\ifplaintex 
\else
\fi

\def\gtm{{\mathsurround=0pt\it $\cal G\mskip-2mu$eometry \&\ 
$\cal T\!\!$opology $\cal M\mskip-1mu$onographs}}    

\def\gtp{{\mathsurround=0pt\it $\cal G\mskip-2mu$eometry \&\ 
$\cal T\!\!$opology $\cal P\!$ublications}}  

\def\recd{{\small Received:\qua\receiveddate\ifx\reviseddate\relax
\else\qquad Revised:\qua\reviseddate\fi\par}} 

\def\volumenumber#1{\def\thevolumenumber{#1}}
\def\volumeyear#1{\def\thevolumeyear{#1}}
\def\volumename#1{\def\thevolumename{#1}}
\def\papernumber#1{\def\thepapernumber{#1}}
\def\pagenumbers#1#2{\def\startpage{#1}\def\finishpage{#2}}
\def\published#1{\def\publishdate{#1}}
\def\received#1{\def\receiveddate{#1}}
\def\revised#1{\def\reviseddate{#1}}
\def\accepted#1{\def\accepteddate{#1}}

\long\def\asciiabstract#1{\long\def\theasciiabstract{#1}}

\let\\\par
\let\thevolumenumber\relax\let\thepapernumber\relax
\let\thevolumeyear\relax\let\startpage\relax
\let\finishpage\relax\let\publishdate\relax\let\receiveddate\relax
\let\reviseddate\relax\let\accepteddate\relax\let\theasciititle\relax
\let\theasciiauthors\relax
\let\theasciiabstract\relax

\let\theerratum\relax\let\theasciiemail\relax
\let\theshortauthors\relax\let\theshorttitle\relax

\def\startpage{1}\def\finishpage{15}\def\thepapernumber{77}

\volumenumber{2}
\volumename{Proceedings of the Kirbyfest}
\volumeyear{1999}

\long\def\maketitlep{   

\count0=\startpage

\gtm\nl        
{\small Volume \thevolumenumber: \thevolumename\nl 
\ifx\theerratum\relax\else Erratum \erratumnumber\nl\fi
Pages \startpage--\finishpage\nl}

\vglue 0.1truein   

{\parskip=0pt\leftskip 0pt plus 1fil\def\\{\par\smallskip}{\ifplaintex\large
\else\Large\fi\bf\thetitle}\par\medskip}   
\vglue 0.05truein 

{\parskip=0pt\leftskip 0pt plus 1fil\def\\{\par}{\sc\theauthors}
\par\medskip}%
 
\vglue 0.03truein 

{\small\leftskip 25pt\rightskip 25pt{\bf Abstract}\stdspace\theabstract

{\bf AMS Classification}\stdspace\theprimaryclass
\ifx\thesecondaryclass\relax\else; \thesecondaryclass\fi\par
{\bf Keywords}\stdspace \thekeywords\par}\vglue 7pt

}   

\font\phead=cmsl9 scaled 950
\font=cmsl9 scaled 1050
\font\pnum=cmbx10 scaled 913
\font\lnum=cmbx10 
\font\pfoot=cmsl9 scaled 950
\font=cmsl9 scaled 1050
\ifplaintex
\headline{\vbox to 0pt{\vskip -4.5mm\line{\small\phead\ifnum
\count0=\startpage ISSN 1464-8997 (on line)
1464-8989 (printed) \hfill {\pnum\folio}\else\ifodd\count0\def\\{ }%
\ifx\theshorttitle\relax\thetitle\else\theshorttitle\fi\hfill{\pnum\folio}
\else\def\\{ and }{\pnum\folio}\hfill\ifx\theshortauthors\relax\theauthors
\else\theshortauthors\fi\fi\fi}\vss}}
\footline{\vbox to 0pt{\vglue 0mm\line{\small\pfoot\ifnum\count0=\startpage
Published \publishdate:\qua\copyright\ \gtp\hfill\else
\gtm, Volume \thevolumenumber\ (\thevolumeyear)\hfill\fi}\vss
}}
\else
\makeatletter
\def\@oddhead{{\small\lhead\ifnum\count0=\startpage ISSN 1464-8997 (on line)
1464-8989 (printed) \hfill {\lnum\number\count0}\else\ifodd\count0
\def\\{ }\ifx\theshorttitle\relax \thetitle \else\theshorttitle\fi\hfill
{\lnum\number\count0}\else\def\\{ and }{\lnum\number\count0}
\hfill\ifx\theshortauthors\relax 
\theauthors\else\theshortauthors\fi\fi\fi}}\def\@evenhead{@oddhead}
\def\@oddfoot{\small\lfoot\ifnum\count0=\startpage Published \publishdate:\qua\copyright\ \gtp\hfill\else
\gtm, Volume \thevolumenumber\ (\thevolumeyear)\hfill\fi}
\def\@evenfoot{@oddfoot}
\makeatother
\fi

\let\maketitlepage\maketitlep
\let\makeshorttitle\maketitlepage
\let\maketitle\maketitlepage

\newwrite\gtoutfile
\long\gdef\makeheadfile{  
{\def\\{, }\def\s{ }
\immediate\openout\gtoutfile head.xxx
\immediate\write\gtoutfile{To: math@arxiv.org}
\immediate\write\gtoutfile{Subject: put OR rep NNNNN:ppppp}
\immediate\write\gtoutfile{--text follows this line--}
\immediate\write\gtoutfile{Proxy-for: \ifx\theasciiauthors\relax
\theauthors\else\theasciiauthors\fi\s<\ifx\theasciiemail\relax\theemail\else\theasciiemail\fi>}
\immediate\write\gtoutfile{\noexpand\\}
\immediate\write\gtoutfile{Authors: \ifx\theasciiauthors\relax
\theauthors\else\theasciiauthors\fi}
{\def\\{ }\immediate\write\gtoutfile{Title: \ifx\theasciititle\relax
\thetitle\else\theasciititle\fi}}
\immediate\write\gtoutfile{Subj-class: GT or SG, GR etc}
\immediate\write\gtoutfile{MSC-class: \theprimaryclass\ifx\thesecondaryclass\relax\else, \thesecondaryclass\fi}
\immediate\write\gtoutfile{Journal-ref: Algebr. Geom. Topol. \thevolumenumber\s
(\thevolumeyear) \startpage-\finishpage}
\immediate\write\gtoutfile{Comments: Published by Geometry and Topology Monographs at}
\immediate\write\gtoutfile{\s\s\s  http://www.maths.warwick.ac.uk/gt/GTMon\thevolumenumber/paper\thepapernumber.abs.html}
\immediate\write\gtoutfile{\noexpand\\}
\immediate\write\gtoutfile{}
\ifx\theasciiabstract\relax
\immediate\write\gtoutfile{\theabstract}\else
\immediate\write\gtoutfile{\theasciiabstract}\fi
\immediate\write\gtoutfile{}
\immediate\write\gtoutfile{\noexpand\\}
\immediate\write\gtoutfile{}
\immediate\closeout\gtoutfile}}  

\def\maketitlepage{\maketitlep\makeheadfile}
\let\makeshorttitle\maketitlepage
\let\maketitle\maketitlepage

\volumenumber{4}
\volumename{Invariants of knots and 3-manifolds (Kyoto 2001)}
\volumeyear{2002}
\papernumber{1}
\pagenumbers{1}{12}
\received{26 November 2001}
\revised{17 February 2002}
\accepted{22 July 2002}
\published{19 September 2002}

\usepackage{dbnsymb,graphicx,curves,epic,eepic,floatflt,picins,amsmath,amssymb}

\theoremstyle{plain}

\newtheorem{theorem}{Theorem}
\newtheorem{proposition}{Proposition}[section]

\theoremstyle{definition}
\newtheorem{definition}[proposition]{Definition}

\theoremstyle{remark}

\newtheorem{remark}[proposition]{Remark}

\newcommand{\mathmode}[1]{$#1$}
\newlength{\standardunitlength}
\newcommand{\im}{\operatorname{im}}

\newlength{\globalparindent}
\def\bbQ{{\mathbb Q}}
\def\bbR{{\mathbb R}}
\def\calA{{\mathcal A}}
\def\calD{{\mathcal D}}
\def\calK{{\mathcal K}}

\def\FI{\mathit{FI}}
\def\FourT{\mathit{4T}}
\def\lin{\mathit{lin}}
\def\GFI{\mathit{GFI}}
\def\GFourT{\mathit{G4T}}
\def\TFI{\mathit{TFI}}
\def\TFourT{\mathit{T4T}}

\begin{document}
\newdimen\captionwidth\captionwidth=\hsize

\title{Bracelets and the Goussarov Filtration\\of the Space of Knots}

\author{Dror Bar-Natan}
\address{Department of Mathematics, University of Toronto\\Toronto 
Ontario M5S 3G3, Canada}
\email{drorbn@math.toronto.edu}
\url{http://www.math.toronto.edu/\char'176drorbn}


\begin{abstract}
Following Goussarov's paper ``Interdependent Modifications of Links and
Invariants of Finite Degree''~\cite{Goussarov:InterdependentModifications}
we describe an alternative finite type theory of knots. While (as shown
by Goussarov) the alternative theory turns out to be equivalent to the
standard one, it nevertheless has its own share of intrinsic beauty.
\end{abstract}

\asciiabstract{Following Goussarov's paper `Interdependent
Modifications of Links and Invariants of Finite Degree' [Topology 37
(1998) 595--602] we describe an alternative finite type theory of
knots. While (as shown by Goussarov) the alternative theory turns out
to be equivalent to the standard one, it nevertheless has its own
share of intrinsic beauty.}

\maketitle
\cl{\small\em In Memory of Mikhail Nikolaevitch Goussarov}

\section{Introduction} There is a well known notion of Vassiliev finite
type invariants of knots (see e.g.~\cite{Bar-Natan:OnVassiliev}). A knot
invariant $I$ is called ``Vassiliev of type $n$'' if, like a polynomial of
degree $n$, its higher than $n$th iterated differences (``derivatives'')
vanish. That is, one picks a knot $K$ and (say) some number $m>n$ of
crossings and then looks at the alternating sum of the values of $I$
evaluated on the $2^m$ knots obtained from $K$ by flipping the crossings
in some subset of the $m$ chosen crossings (with signs determined by the
parity of the number of crossings flipped). If this sum vanishes for all
$K$ and all choices of $m>n$ crossings, then $I$ is of Vassiliev type $n$.

A different way of saying this is to say that we look at $K$ and at some
number $m$ of possible simple modifications to $K$ (of the form
$\overcrossing\to\undercrossing$ or $\undercrossing\to\overcrossing$) which
can (but don't need to) be performed simultaneously. We then look at
iterated differences of values of $I$ evaluated on $K$ with just some of
these modification applied, and if this vanishes whenever $m>n$, then $I$
is of Vassiliev type $n$.

But why restrict to just ``simple modification''? Goussarov's novel
idea in his paper `Interdependent Modifications of Links and Invariants
of Finite Degree''~\cite{Goussarov:InterdependentModifications} was to
allow arbitrary modifications to $K$. That is, we pick some number $m$ of
intervals along $K$ and allow them to make completely arbitrary detours,
provided none of the original paths and none of the re-routed paths ever
intersect. We can then form the same sort of alternating sum of values of
a knot invariant $I$, and make a similar definition of ``Goussarov type
$n$'', if this alternating sum vanishes whenever the number of detours $m$
is bigger than $n$. (We will repeat this definition in more precise terms
in Section~\ref{sec:Goussarov}).

\begin{figure}[ht!]
  \begin{center}
    \includegraphics[width=2in]{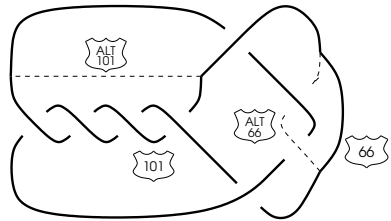}
  \end{center}
  \caption{$6_1$ and two detours}\label{fig:Rerouting}
  \end{figure}

An example appears in Figure~\ref{fig:Rerouting}; if we travel the
main road, it is the knot $6_1$. If we choose route \textsc{alt~66}
over route \textsc{66}, the knot becomes the more complicated
$8_3$. If we choose route \textsc{alt~101} over route \textsc{101} we
get the unknot no matter which choice we make in the east. Thus the
alternating sum corresponding to this knot and this choice of detours
is $I(6_1)-I(8_3)-I(0)+I(0)$.

Goussarov's theorem says that the two notions of finite type invariants
agree up to some renumbering:

\begin{theorem}[Goussarov~\cite{Goussarov:InterdependentModifications}]
\label{thm:main}
Any Vassiliev type $n$ invariant is a Goussarov type $2n$ invariant and any
Goussarov type $2n$ or $2n+1$ invariant is a Vassiliev type $n$ invariant.
\end{theorem}

\begin{floatingfigure}[l]{2in}
\begin{center}
  \includegraphics[width=1.8in]{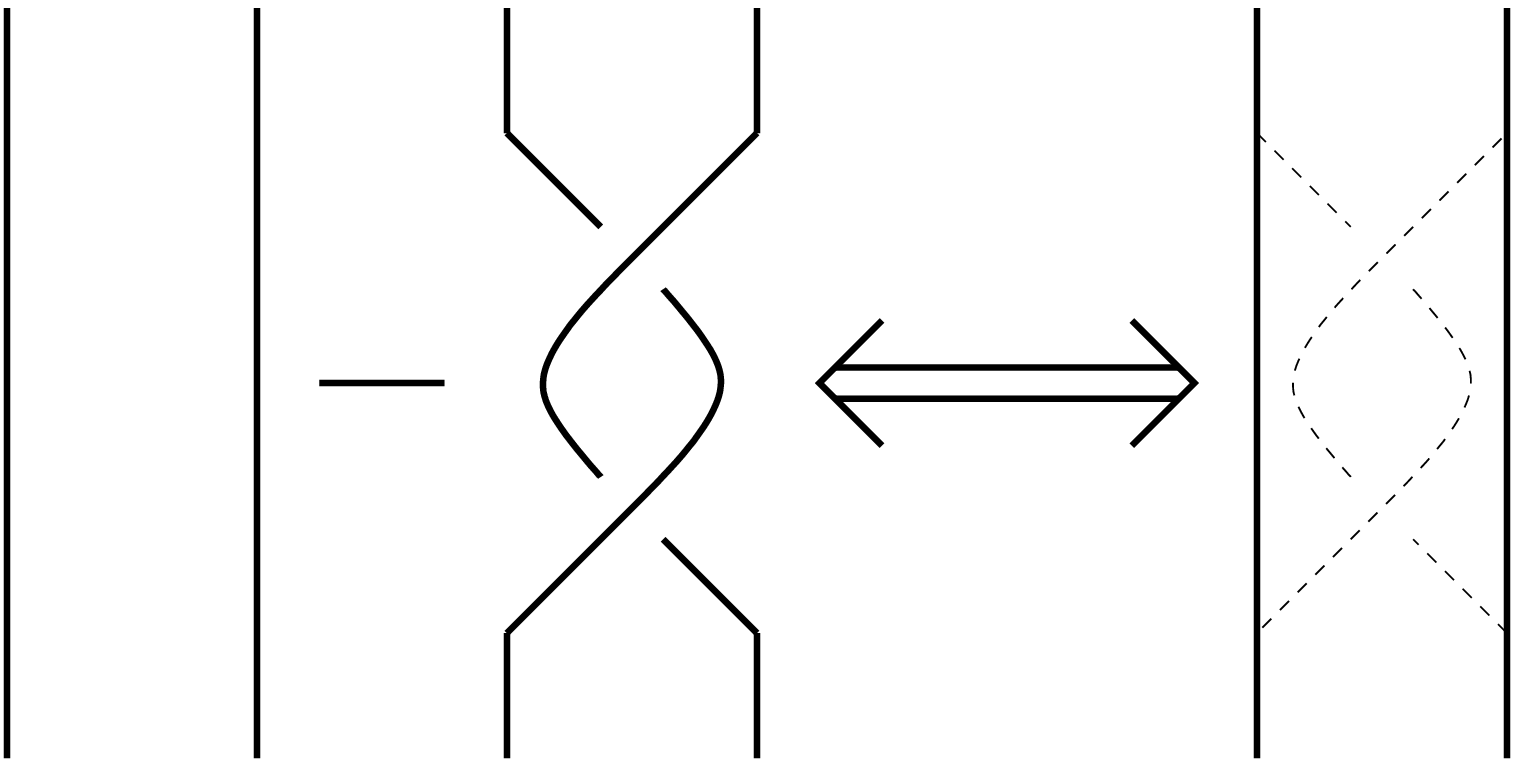}
\end{center}
\end{floatingfigure}
The key to the understanding of this theorem is the figure on the left,
which indicates that a single Vassiliev style crossing change (left part
of the figure) can be achieved using two Goussarov style detour moves
(right part of the figure). Indeed, if none or just one of the detours
is taken, the knot-part displayed remains unbraided, and only if both
detours are taken do we get braiding. This too will be made precise later
in this paper.

Our paper is only partially about proving Theorem~\ref{thm:main}. The
theorem says that the two notions of finite type invariants are
equivalent. Thus if we start from the Goussarov notion and study it
along the same lines as the standard study of the Vassiliev notion, we
must meet the same objects: chord diagrams, $\FourT$ relations, etc.,
even if we pretend to know nothing about Vassiliev finite type
invariants and about Theorem~\ref{thm:main}. Hence our plan is to carry
out an independent study of the Goussarov theory with the hope that we
encounter some familiar objects as we go. This we do in
Section~\ref{sec:Goussarov} which to our taste is the most elegant part
of this paper. Before that, in Section~\ref{sec:Vassiliev}, we quickly
review the basics of the Vassiliev theory. This review is not a
prerequisite for the study of the Goussarov theory (or else we would be
defeating our own purpose), and we embark upon it merely for the
purpose of comparison and to establish what we mean by the word
``study''. Finally, in Section~\ref{sec:Equivalence} we use some of the
results of Section~\ref{sec:Goussarov} to give an easy proof of
Theorem~\ref{thm:main}.

A different easy proof of Theorem~\ref{thm:main} is in
Conant's~\cite{Conant:OnGoussarov}.

\medskip

{\bf Acknowledgements}\qua This article is the written form
of a lecture given by the author at the Goussarov Day at RIMS, Kyoto,
September 25, 2001. It follows research done jointly with my student
Haggai Scolnicov in 1998/99 and conversations with M.~Hutchings in 1997.
I wish to thank Shlomi Kotler for pointing out an error in a previous
version of Figure~\ref{fig:G4TFI} and A.~Referee for further comments.

This paper is also available electronically at:\qua {\tt arXiv:math.GT/0111267}
\nl and at:\qua {\tt
  http://www.math.toronto.edu/\char'176drorbn/papers/Bracelets/}

\section{A quick review of the Vassiliev finite type theory}
\label{sec:Vassiliev} The purpose of this section is to recall how chord
diagrams and the $\FourT$ relations arise in the Vassiliev theory of finite
type invariants.

Let $\calK^V_n$ denote the space of all formal linear combinations
of $n$-singular knots, knots with $n$ ``double points''
that locally look like $\doublepoint$, modulo the benign
``differentiability relations'' which will be described shortly. Let
$\delta^V=\delta^V_{n+1}:\calK^V_{n+1}\to\calK^V_n$ be the linear map defined on
a singular knot $K$ by picking one of the double points $\doublepoint$
in $K$ and then mapping $K$ to the difference of the knots obtained by
resolving $\doublepoint$ to and overcrossing $\overcrossing$ and to an
undercrossing $\undercrossing$:
\[ \delta^V:\doublepoint\mapsto\overcrossing-\undercrossing. \]
As it stands, $\delta^V$ is not well defined because it may depend on the
choice of the double point to be resolved. We fix this by dividing
$\calK^V_n$ by differentiability relations, which are exactly the minimal
relations required in order to make $\delta^V$ well defined. In figures, the
differentiability relations are the relations
\[
  \overcrossing\doublepoint-\undercrossing\doublepoint
  = \doublepoint\overcrossing-\doublepoint\undercrossing.
\]
(As usual in knot theory, this equation represents the whole family of
relations obtained from the figures drawn by completing them to knots in
all possible ways, but where all the ``picturelets'' (like
$\overcrossing\doublepoint$ and $\doublepoint\undercrossing$) are completed
in the same manner).

We denote the adjoint of $\delta^V$ by $\partial_V$
and call it ``the derivative''. It is a map
$\partial_V:(\calK^V_n)^\star\to(\calK^V_{n+1})^\star$. The name
``derivative'' is justified by the fact that $(\partial_V I)(K)$ for
some $I\in(\calK^V_n)^\star$ and some generator $K\in\calK^V_{n+1}$ is by
definition the difference of the values of $I$ on two ``neighboring''
$n$-singular knots, in harmony with the usual definition of derivatives
for functions on $\bbR^d$.

\begin{definition} A knot invariant $I$ (equivalently, a linear
functional on $\calK=\calK^V_0$) is of Vassiliev type $n$ if
its $(n+1)$-st (Vassiliev style) derivative vanishes, that is, if
$(\partial_V)^{n+1}I\equiv 0$.  (This definition is the analog of one
of the standard definitions of polynomials on $\bbR^d$).
\end{definition}

When thinking about finite type invariants, it is convenient to have in
mind the following ladders of spaces and their duals, printed here with the
names of some specific elements that we will use later:
\begin{equation} \label{eq:Vladders}
    \def\neg{\hspace{-12pt}}
    \begin{array}{cccccccccr}
      \ldots\longrightarrow\neg &
        \calK^V_{n+1} &
        \neg\stackrel{\delta^V}{\longrightarrow}\neg &
        \calK^V_n &
        \neg\stackrel{\delta^V}{\longrightarrow} &
        \calK^V_{n-1} &
        \longrightarrow\ldots &
        \stackrel{\delta^V}{\longrightarrow} &
        \calK^V_0=\calK \\
      &&&&&&&& \\
      \ldots\longleftarrow\neg &
        (\calK^V_{n+1})^\star &
        \neg\stackrel{\partial_V}{\longleftarrow}\neg &
        (\calK^V_n)^\star &
        \neg\stackrel{\partial_V}{\longleftarrow} &
        (\calK^V_{n-1})^\star &
        \longleftarrow\ldots &
        \stackrel{\partial_V}{\longleftarrow} &
        (\calK^V_0)^\star=\calK^\star \\
      & \inup\rule{0pt}{16pt}     & & \inup &
        & & & & \inup\  \\
      & \partial_V^{n+1}I\equiv 0   & & \partial_V^nI=W &
        & & & & I\
    \end{array}
  \end{equation}
We often study invariants of type $n$ by
studying their $n$th derivatives. Clearly, if $I$ is of type $n$ and
$W=\partial_V^nI$, then $\partial_V W=0$ (``$W$ is a constant''). Glancing
at~\eqref{eq:Vladders}, we see that $W$ descends to a linear functional,
also called $W$, on $\calK^V_n/\delta^V\calK^V_{n+1}$. The latter space is a
familiar entity:

\begin{proposition} \label{prop:Vsymbols}
The space $\calK^V_n/\delta^V\calK^V_{n+1}$ is canonically isomorphic to the
space $\calD^V_n$ of $n$-chord diagrams, defined below. \qed
\end{proposition}

\begin{definition} An $n$-chord diagram is a choice of $n$ pairs of
distinct points on an oriented circle, considered up to orientation
preserving homeomorphisms of the circle. Usually an $n$-chord diagram is
simply drawn as a circle with $n$ chords (whose ends are the $n$ pairs).
The space $\calD^V_n$ is the space of all formal linear combinations of
$n$-chord diagrams. As an example, a basis for $\calD^V_3$ is
$\{
  \parbox{4mm}{\setlength{\unitlength}{2mm} \begin{picture}(2,2)(-1,-1)
  \put(0,0){\bigcircle{2}} \curve(0.17365,-0.98481,1.,0.)
  \curve(0.76604,0.64279,-0.5,0.86603)
  \curve(-0.93969,0.34202,-0.5,-0.86603)
  \end{picture}}, \parbox{4mm}{\setlength{\unitlength}{2mm}
  \begin{picture}(2,2)(-1,-1) \put(0,0){\bigcircle{2}}
  \curve(-0.22252,-0.97493,1.,0.) \curve(0.62349,0.78183,-0.90097,0.43388)
  \curve(-0.90097,-0.43388,0.62349,-0.78183) \end{picture}},
  \parbox{4mm}{\setlength{\unitlength}{2mm} \begin{picture}(2,2)(-1,-1)
  \put(0,0){\bigcircle{2}} \curve(-0.70711,-0.70711,0.70711,-0.70711)
  \curve(0.70711,0.70711,-0.70711,0.70711) \curve(-1.,0.,1.,0.)
  \end{picture}}, \parbox{4mm}{\setlength{\unitlength}{2mm}
  \begin{picture}(2,2)(-1,-1) \put(0,0){\bigcircle{2}}
  \curve(-0.5,0.86603,0.5,-0.86603) \curve(-0.5,-0.86603,1.,0.)
  \curve(0.5,0.86603,-1.,0.) \end{picture}},
  \parbox{4mm}{\setlength{\unitlength}{2mm} \begin{picture}(2,2)(-1,-1)
  \put(0,0){\bigcircle{2}} \curve(-0.5,0.86603,0.5,-0.86603)
  \curve(0.5,0.86603,-0.5,-0.86603) \curve(-1.,0.,1.,0.) \end{picture}}
\}$.
\end{definition}

Next, we wish to find conditions that a ``potential top derivative''
has to satisfy in order to actually be a top derivative. More precisely,
we wish to find conditions that a functional  $W\in(\calD^V_n)^\star$
has to satisfy in order to be $\partial_V^nI$ for some invariant $I$. A
first condition is that $W$ must be ``integrable once''; namely,
there has to be some $W^1\in(\calK^V_{n-1})^V$ with $W=\partial_V
W^1$. Another quick glance at~\eqref{eq:Vladders}, and we see that $W$
is integrable once iff it vanishes on $\ker\delta^V$, which is the same as
requiring that $W$ descends to $\calA^V_n := \calD^V_n/\pi(\ker\delta^V)
= \calK^V_n/(\im\delta^V+\ker\delta^V)$ ($\pi$ is the projection
$\calK^V_n\to\calD^V_n=\calK^V_n/\delta^V\calK^V_{n+1}$, and there
should be no confusion regarding the identities of the $\delta^V$'s
involved). Often elements of $(\calA^V_n)^\star$ are referred to as
``weight systems''. A more accurate name would be ``once-integrable
weight systems''.

We see that it is necessary to understand $\ker\delta^V$. In
Figure~\ref{fig:T4T} we show a family of members of $\ker\delta^V$, the
``Topological 4-Term'' ($\TFourT$) relations.  Figure~\ref{fig:Lasso} explains
how they arise from ``lassoing a singular point''. Figure~\ref{fig:TFI}
shows another family of members of $\ker\delta^V$. The following theorem
says that this is all:

\begin{figure}[ht!]
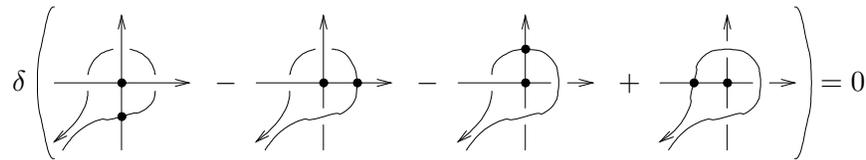

\[ 
  \setlength{\unitlength}{0.7\standardunitlength}
  \begin{array}{c}
    {\input figs/T4T.tex } 
  \end{array}
 \]
\caption{
  A Topological 4-Term ($\TFourT$) relation.
  Each of the four graphics in the
  picture represents a part of an $n$-singular knot (so there are $n-2$
  additional singular points not shown), and, as usual in knot theory, the
  4 singular knots in the equation are the same outside the region shown.
} \label{fig:T4T}
\end{figure}

\begin{figure}[ht!]
\[ 
  \setlength{\unitlength}{0.7\standardunitlength}
  \begin{array}{c}
    {\begingroup\makeatletter\ifx\SetFigFont\undefined%
\gdef\SetFigFont#1#2#3#4#5{%
  \reset@font\fontsize{#1}{#2pt}%
  \fontfamily{#3}\fontseries{#4}\fontshape{#5}%
  \selectfont}%
\fi\endgroup%
{\renewcommand{\dashlinestretch}{30}
\begin{picture}(7524,1239)(0,-10)
\put(2712,612){\blacken\ellipse{70}{70}}
\put(2712,612){\ellipse{70}{70}}
\path(2712,12)(2712,1212)
\path(2742.000,1092.000)(2712.000,1212.000)(2682.000,1092.000)
\path(2112,612)(3312,612)
\path(3192.000,582.000)(3312.000,612.000)(3192.000,642.000)
\path(3012,687)(3011,690)(3010,696)
	(3008,705)(3005,717)(3000,732)
	(2995,747)(2988,763)(2980,781)
	(2969,799)(2954,818)(2937,837)
	(2918,854)(2899,869)(2881,880)
	(2863,888)(2847,895)(2832,900)
	(2817,905)(2805,908)(2796,910)
	(2790,911)(2787,912)
\path(2637,912)(2634,911)(2628,910)
	(2619,908)(2607,905)(2592,900)
	(2577,895)(2561,888)(2543,880)
	(2525,869)(2506,854)(2487,837)
	(2470,818)(2455,799)(2444,781)
	(2436,763)(2429,747)(2424,732)
	(2419,717)(2416,705)(2414,696)
	(2413,690)(2412,687)
\path(2412,537)(2412,535)(2412,530)
	(2411,522)(2411,511)(2409,497)
	(2407,481)(2404,463)(2400,443)
	(2393,422)(2384,399)(2372,373)
	(2357,344)(2337,312)(2317,284)
	(2297,258)(2277,235)(2257,213)
	(2238,194)(2219,176)(2200,159)
	(2182,143)(2165,129)(2150,116)
	(2136,106)(2112,87)
\path(2187.464,185.006)(2112.000,87.000)(2224.707,137.963)
\path(2787,312)(2790,313)(2796,314)
	(2805,316)(2817,319)(2832,324)
	(2847,329)(2863,336)(2881,344)
	(2899,355)(2918,370)(2937,387)
	(2954,406)(2969,425)(2980,443)
	(2988,461)(2995,477)(3000,492)
	(3005,507)(3008,519)(3010,528)
	(3011,534)(3012,537)
\path(2637,312)(2635,312)(2630,312)
	(2622,311)(2611,311)(2597,309)
	(2581,307)(2563,304)(2543,300)
	(2522,293)(2499,284)(2473,272)
	(2444,257)(2412,237)(2384,217)
	(2358,197)(2335,177)(2313,157)
	(2294,138)(2276,119)(2259,100)
	(2243,82)(2229,65)(2216,50)
	(2206,36)(2197,26)(2192,18)
	(2189,14)(2187,12)
\put(612,612){\blacken\ellipse{70}{70}}
\put(612,612){\ellipse{70}{70}}
\path(12,612)(1212,612)
\path(1092.000,582.000)(1212.000,612.000)(1092.000,642.000)
\path(612,12)(612,1212)
\path(642.000,1092.000)(612.000,1212.000)(582.000,1092.000)
\path(87,12)(88,13)(91,17)
	(95,22)(101,29)(108,39)
	(117,50)(128,61)(140,74)
	(154,88)(170,104)(189,121)
	(211,141)(237,162)(261,180)
	(283,197)(303,210)(319,220)
	(333,227)(344,232)(353,235)
	(362,237)(370,239)(379,242)
	(389,247)(401,254)(415,264)
	(430,277)(447,294)(462,312)
	(473,328)(481,344)(487,357)
	(491,368)(494,377)(495,384)
	(495,390)(494,395)(493,400)
	(492,404)(490,408)(488,413)
	(486,419)(483,426)(480,434)
	(475,443)(469,453)(462,462)
	(453,469)(443,475)(434,480)
	(426,483)(419,486)(413,488)
	(408,490)(404,492)(399,493)
	(395,494)(390,495)(384,495)
	(377,494)(368,491)(357,487)
	(344,481)(328,473)(312,462)
	(294,447)(277,430)(264,415)
	(254,401)(247,389)(242,379)
	(239,370)(237,362)(235,353)
	(232,344)(227,333)(220,319)
	(210,303)(197,283)(180,261)
	(162,237)(141,211)(121,189)
	(104,170)(88,154)(74,140)
	(61,128)(49,117)(39,108)
	(29,101)(12,87)
\path(85.560,186.443)(12.000,87.000)(123.703,140.127)
\put(4812,612){\blacken\ellipse{70}{70}}
\put(4812,612){\ellipse{70}{70}}
\path(4587,612)(5037,612)
\path(4812,387)(4812,837)
\path(4437,612)(4212,612)
\path(4812,987)(4812,1212)
\path(4842.000,1092.000)(4812.000,1212.000)(4782.000,1092.000)
\path(5187,612)(5412,612)
\path(5292.000,582.000)(5412.000,612.000)(5292.000,642.000)
\path(4812,237)(4812,12)
\path(4287,12)(4288,14)(4291,18)
	(4296,26)(4303,36)(4312,50)
	(4323,65)(4336,82)(4351,100)
	(4366,119)(4384,138)(4404,157)
	(4426,177)(4451,197)(4480,217)
	(4512,237)(4545,254)(4575,268)
	(4600,278)(4620,284)(4635,288)
	(4647,289)(4655,288)(4662,287)
	(4669,285)(4677,284)(4689,284)
	(4704,286)(4724,289)(4749,295)
	(4779,303)(4812,312)(4848,322)
	(4878,331)(4902,336)(4919,339)
	(4930,340)(4938,339)(4943,337)
	(4948,335)(4955,335)(4964,337)
	(4977,343)(4994,353)(5015,367)
	(5037,387)(5054,406)(5069,425)
	(5080,442)(5087,456)(5093,466)
	(5096,475)(5098,481)(5100,487)
	(5100,493)(5101,500)(5103,510)
	(5104,522)(5107,539)(5109,559)
	(5111,584)(5112,612)(5111,640)
	(5109,665)(5107,685)(5104,702)
	(5103,714)(5101,724)(5100,731)
	(5100,737)(5098,743)(5096,749)
	(5093,758)(5087,768)(5080,782)
	(5069,799)(5054,818)(5037,837)
	(5018,854)(4999,869)(4982,880)
	(4968,887)(4958,893)(4949,896)
	(4943,898)(4937,899)(4931,900)
	(4924,901)(4914,903)(4902,904)
	(4885,907)(4865,909)(4840,911)
	(4812,912)(4784,911)(4759,909)
	(4739,907)(4722,904)(4710,903)
	(4700,901)(4693,900)(4687,899)
	(4681,898)(4675,896)(4666,893)
	(4656,887)(4642,880)(4625,869)
	(4606,854)(4587,837)(4567,815)
	(4553,794)(4543,777)(4537,764)
	(4535,755)(4535,748)(4537,743)
	(4539,738)(4540,730)(4539,719)
	(4536,702)(4531,678)(4522,648)
	(4512,612)(4503,579)(4495,549)
	(4489,524)(4486,504)(4484,489)
	(4484,477)(4485,469)(4487,462)
	(4488,455)(4489,447)(4488,435)
	(4484,420)(4478,400)(4468,375)
	(4454,345)(4437,312)(4417,280)
	(4397,251)(4377,226)(4357,204)
	(4338,184)(4319,166)(4300,151)
	(4282,136)(4265,123)(4250,112)
	(4236,103)(4212,87)
\path(4295.205,178.526)(4212.000,87.000)(4328.487,128.603)
\put(6912,612){\blacken\ellipse{70}{70}}
\put(6912,612){\ellipse{70}{70}}
\path(6312,612)(7512,612)
\path(7392.000,582.000)(7512.000,612.000)(7392.000,642.000)
\path(6912,12)(6912,1212)
\path(6942.000,1092.000)(6912.000,1212.000)(6882.000,1092.000)
\path(6387,12)(6388,13)(6391,17)
	(6395,22)(6401,29)(6408,39)
	(6417,50)(6428,61)(6440,74)
	(6454,88)(6470,104)(6489,121)
	(6511,141)(6537,162)(6561,180)
	(6583,197)(6603,210)(6619,220)
	(6633,227)(6644,232)(6653,235)
	(6662,237)(6670,239)(6679,242)
	(6689,247)(6701,254)(6715,264)
	(6730,277)(6747,294)(6762,312)
	(6773,328)(6781,344)(6787,357)
	(6791,368)(6794,377)(6795,384)
	(6795,390)(6794,395)(6793,400)
	(6792,404)(6790,408)(6788,413)
	(6786,419)(6783,426)(6780,434)
	(6775,443)(6769,453)(6762,462)
	(6753,469)(6743,475)(6734,480)
	(6726,483)(6719,486)(6713,488)
	(6708,490)(6704,492)(6699,493)
	(6695,494)(6690,495)(6684,495)
	(6677,494)(6668,491)(6657,487)
	(6644,481)(6628,473)(6612,462)
	(6594,447)(6577,430)(6564,415)
	(6554,401)(6547,389)(6542,379)
	(6539,370)(6537,362)(6535,353)
	(6532,344)(6527,333)(6520,319)
	(6510,303)(6497,283)(6480,261)
	(6462,237)(6441,211)(6421,189)
	(6404,170)(6388,154)(6374,140)
	(6361,128)(6349,117)(6339,108)
	(6329,101)(6312,87)
\path(6385.560,186.443)(6312.000,87.000)(6423.703,140.127)
\path(1512,612)(1812,612)
\path(1692.000,582.000)(1812.000,612.000)(1692.000,642.000)
\path(3612,612)(3912,612)
\path(3792.000,582.000)(3912.000,612.000)(3792.000,642.000)
\path(5712,612)(6012,612)
\path(5892.000,582.000)(6012.000,612.000)(5892.000,642.000)
\end{picture}
} } 
  \end{array}
 \]
\caption{
  Lassoing a singular point: Each of the graphics represents an
  $(n-1)$-singular knot, but only one of the singularities is explicitly
  displayed.  Start from the left-most graphic, pull the ``lasso'' under
  the displayed singular point, ``lasso'' the singular point by crossing
  each of the four arcs emanating from it one at a time, and pull the
  lasso back out, returning to the initial position. Each time an arc is
  crossed, the difference between ``before'' and ``after'' is the
  $\delta^V$ applied to an $n$-singular knot (up to signs). The four
  $n$-singular knot thus obtained are the ones making the Topological
  4-Term relation, and $\delta^V$ applied to their signed sum is the
  difference between the first and the last $(n-1)$-singular knot shown
  in this figure; namely, it is $0$.
} \label{fig:Lasso}
\end{figure}
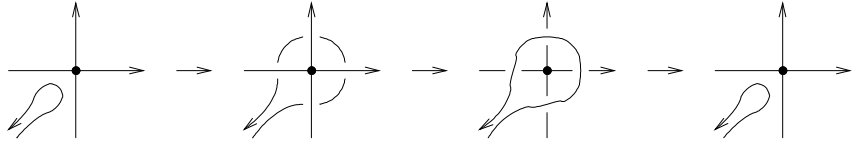

\begin{figure}[ht!]
\[ 
  \setlength{\unitlength}{0.4\standardunitlength}
  \begin{array}{c}
    {\begingroup\makeatletter\ifx\SetFigFont\undefined%
\gdef\SetFigFont#1#2#3#4#5{%
  \reset@font\fontsize{#1}{#2pt}%
  \fontfamily{#3}\fontseries{#4}\fontshape{#5}%
  \selectfont}%
\fi\endgroup%
{\renewcommand{\dashlinestretch}{30}
\begin{picture}(8267,1456)(0,-10)
\path(3225,158)(3226,159)(3227,161)
	(3230,165)(3235,172)(3241,181)
	(3250,192)(3260,207)(3273,223)
	(3287,242)(3303,264)(3321,286)
	(3340,311)(3361,336)(3383,362)
	(3406,389)(3430,417)(3456,445)
	(3483,474)(3512,503)(3542,532)
	(3575,563)(3610,594)(3648,626)
	(3688,659)(3732,692)(3777,725)
	(3825,758)(3873,789)(3920,818)
	(3965,845)(4006,868)(4042,889)
	(4074,907)(4102,923)(4126,936)
	(4147,947)(4164,956)(4178,964)
	(4191,971)(4202,977)(4213,983)
	(4223,988)(4234,993)(4246,998)
	(4260,1004)(4277,1009)(4296,1015)
	(4319,1022)(4345,1028)(4376,1035)
	(4410,1042)(4448,1048)(4489,1053)
	(4532,1057)(4575,1058)(4620,1056)
	(4662,1052)(4701,1046)(4735,1038)
	(4765,1029)(4792,1020)(4814,1011)
	(4834,1002)(4851,993)(4865,984)
	(4878,976)(4889,967)(4900,958)
	(4910,949)(4920,939)(4930,929)
	(4940,918)(4951,906)(4962,893)
	(4973,878)(4985,862)(4996,844)
	(5007,824)(5016,803)(5023,781)
	(5025,758)(5022,733)(5015,709)
	(5004,687)(4992,667)(4978,648)
	(4965,632)(4951,617)(4938,605)
	(4925,593)(4913,583)(4900,573)
	(4887,564)(4874,555)(4860,547)
	(4845,538)(4827,529)(4807,519)
	(4785,509)(4759,499)(4729,489)
	(4696,479)(4659,470)(4618,463)
	(4575,458)(4531,457)(4488,459)
	(4445,463)(4405,470)(4367,479)
	(4330,490)(4295,502)(4262,516)
	(4230,530)(4198,545)(4168,561)
	(4139,578)(4111,594)(4085,610)
	(4062,625)(4040,639)(4022,651)
	(4006,661)(3994,670)(3985,676)
	(3980,680)(3976,682)(3975,683)
\path(3750,833)(3748,834)(3745,836)
	(3739,839)(3729,845)(3717,852)
	(3703,861)(3686,871)(3668,882)
	(3648,895)(3628,909)(3607,924)
	(3584,941)(3560,959)(3535,980)
	(3508,1003)(3479,1030)(3450,1058)
	(3422,1087)(3395,1116)(3372,1143)
	(3351,1168)(3333,1192)(3316,1215)
	(3301,1236)(3287,1256)(3274,1276)
	(3263,1294)(3253,1311)(3244,1325)
	(3237,1337)(3231,1347)(3228,1353)
	(3226,1356)(3225,1358)
\path(5625,1358)(5626,1357)(5627,1355)
	(5630,1351)(5635,1344)(5641,1335)
	(5650,1324)(5660,1309)(5673,1293)
	(5687,1274)(5703,1252)(5721,1230)
	(5740,1205)(5761,1180)(5783,1154)
	(5806,1127)(5830,1099)(5856,1071)
	(5883,1042)(5912,1013)(5942,984)
	(5975,953)(6010,922)(6048,890)
	(6088,857)(6132,824)(6177,791)
	(6225,758)(6273,727)(6320,698)
	(6365,671)(6406,648)(6442,627)
	(6474,609)(6502,593)(6526,580)
	(6547,569)(6564,560)(6578,552)
	(6591,545)(6602,539)(6613,533)
	(6623,528)(6634,523)(6646,518)
	(6660,512)(6677,507)(6696,501)
	(6719,494)(6745,488)(6776,481)
	(6810,474)(6848,468)(6889,463)
	(6932,459)(6975,458)(7020,460)
	(7062,464)(7101,470)(7135,478)
	(7165,487)(7192,496)(7214,505)
	(7234,514)(7251,523)(7265,532)
	(7278,540)(7289,549)(7300,558)
	(7310,567)(7320,577)(7330,587)
	(7340,598)(7351,610)(7362,623)
	(7373,638)(7385,654)(7396,672)
	(7407,692)(7416,713)(7423,735)
	(7425,758)(7422,783)(7415,807)
	(7404,829)(7392,849)(7378,868)
	(7365,884)(7351,899)(7338,911)
	(7325,923)(7313,933)(7300,943)
	(7287,952)(7274,961)(7260,969)
	(7245,978)(7227,987)(7207,997)
	(7185,1007)(7159,1017)(7129,1027)
	(7096,1037)(7059,1046)(7018,1053)
	(6975,1058)(6931,1059)(6888,1057)
	(6845,1053)(6805,1046)(6767,1037)
	(6730,1026)(6695,1014)(6662,1000)
	(6630,986)(6598,971)(6568,955)
	(6539,938)(6511,922)(6485,906)
	(6462,891)(6440,877)(6422,865)
	(6406,855)(6394,846)(6385,840)
	(6380,836)(6376,834)(6375,833)
\path(6150,683)(6148,682)(6145,680)
	(6139,677)(6129,671)(6117,664)
	(6103,655)(6086,645)(6068,634)
	(6048,621)(6028,607)(6007,592)
	(5984,575)(5960,557)(5935,536)
	(5908,513)(5879,486)(5850,458)
	(5822,429)(5795,400)(5772,373)
	(5751,348)(5733,324)(5716,301)
	(5701,280)(5687,260)(5674,240)
	(5663,222)(5653,205)(5644,191)
	(5637,179)(5631,169)(5628,163)
	(5626,160)(5625,158)
\path(2700,833)(3075,833)
\path(2700,683)(3075,683)
\path(7650,833)(8025,833)
\path(7650,683)(8025,683)
\put(881.250,758.000){\arc{3187.500}{5.8458}{6.7205}}
\put(2062.500,720.500){\arc{3525.798}{2.7255}{3.5577}}
\put(1125,758){\blacken\ellipse{150}{150}}
\put(1125,758){\ellipse{150}{150}}
\path(5175,758)(5475,758)
\path(525,1358)(526,1357)(527,1355)
	(530,1351)(535,1344)(541,1335)
	(550,1324)(560,1309)(573,1293)
	(587,1274)(603,1252)(621,1230)
	(640,1205)(661,1180)(683,1154)
	(706,1127)(730,1099)(756,1071)
	(783,1042)(812,1013)(842,984)
	(875,953)(910,922)(948,890)
	(988,857)(1032,824)(1077,791)
	(1125,758)(1173,727)(1220,698)
	(1265,671)(1306,648)(1342,627)
	(1374,609)(1402,593)(1426,580)
	(1447,569)(1464,560)(1478,552)
	(1491,545)(1502,539)(1513,533)
	(1523,528)(1534,523)(1546,518)
	(1560,512)(1577,507)(1596,501)
	(1619,494)(1645,488)(1676,481)
	(1710,474)(1748,468)(1789,463)
	(1832,459)(1875,458)(1920,460)
	(1962,464)(2001,470)(2035,478)
	(2065,487)(2092,496)(2114,505)
	(2134,514)(2151,523)(2165,532)
	(2178,540)(2189,549)(2200,558)
	(2210,567)(2220,577)(2230,587)
	(2240,598)(2251,610)(2262,623)
	(2273,638)(2285,654)(2296,672)
	(2307,692)(2316,713)(2323,735)
	(2325,758)(2323,781)(2316,803)
	(2307,824)(2296,844)(2285,862)
	(2273,878)(2262,893)(2251,906)
	(2240,918)(2230,929)(2220,939)
	(2210,949)(2200,958)(2189,967)
	(2178,976)(2165,984)(2151,993)
	(2134,1002)(2114,1011)(2092,1020)
	(2065,1029)(2035,1038)(2001,1046)
	(1962,1052)(1920,1056)(1875,1058)
	(1832,1057)(1789,1053)(1748,1048)
	(1710,1042)(1676,1035)(1645,1028)
	(1619,1022)(1596,1015)(1577,1009)
	(1560,1004)(1546,998)(1534,993)
	(1523,988)(1512,983)(1502,977)
	(1491,971)(1478,964)(1464,956)
	(1447,947)(1426,936)(1402,923)
	(1374,907)(1342,889)(1306,868)
	(1265,845)(1220,818)(1173,789)
	(1125,758)(1077,725)(1032,692)
	(988,659)(948,626)(910,594)
	(875,563)(842,532)(812,503)
	(783,474)(756,445)(730,417)
	(706,389)(683,362)(661,336)
	(640,311)(621,286)(603,264)
	(587,242)(573,223)(560,207)
	(550,192)(541,181)(535,172)
	(530,165)(527,161)(526,159)(525,158)
\put(0,608){\makebox(0,0)[lb]{\smash{{\mathmode{\delta}}}}}
\put(8175,608){\makebox(0,0)[lb]{\smash{{\mathmode{0}}}}}
\end{picture}
} } 
  \end{array}
 \]
\caption{
  A Topological Framing Independence Relation ($\TFI$)
} \label{fig:TFI}
\end{figure}

\begin{theorem}[Stanford~\cite{Stanford:FiniteType}]\label{thm:Stanford} 
The $\TFourT$ relations of Figure~\ref{fig:T4T} and the $\TFI$ relations of
Figure~\ref{fig:TFI} span $\ker\delta^V$. \qed
\end{theorem}

Pushing the $\TFourT$ and the $\TFI$ relations down to the level of chord
diagrams, we get the well-known $\FourT$ and $\FI$ relations, which span
$\pi(\ker\delta^V)$: (see e.g.~\cite{Bar-Natan:OnVassiliev})
\[ \FourT:\quad
  \setlength{\unitlength}{0.5\standardunitlength}
  \begin{array}{c}
    {\begingroup\makeatletter\ifx\SetFigFont\undefined%
\gdef\SetFigFont#1#2#3#4#5{%
  \reset@font\fontsize{#1}{#2pt}%
  \fontfamily{#3}\fontseries{#4}\fontshape{#5}%
  \selectfont}%
\fi\endgroup%
{\renewcommand{\dashlinestretch}{30}
\begin{picture}(6476,1089)(0,-10)
\put(1288,462){\makebox(0,0)[lb]{\smash{{\mathmode{-}}}}}
\put(3088,462){\makebox(0,0)[lb]{\smash{{\mathmode{=}}}}}
\put(4888,462){\makebox(0,0)[lb]{\smash{{\mathmode{-}}}}}
\put(538,537){\ellipse{1060}{1060}}
\put(2338,537){\ellipse{1060}{1060}}
\put(4138,537){\ellipse{1060}{1060}}
\put(5938,537){\ellipse{1060}{1060}}
\path(913,912)(913,162)
\path(2713,912)(2713,162)
\path(4513,912)(4513,162)
\path(6313,912)(6313,162)
\path(13,537)(838,987)
\path(1813,537)(2788,837)
\path(3613,537)(4438,87)
\path(5413,537)(6388,237)
\end{picture}
} } 
  \end{array}
\qquad \FI:\quad
  \setlength{\unitlength}{0.5\standardunitlength}
  \begin{array}{c}
    {%
\begingroup\makeatletter\ifx\SetFigFont\undefined%
\gdef\SetFigFont#1#2#3#4#5{%
  \reset@font\fontsize{#1}{#2pt}%
  \fontfamily{#3}\fontseries{#4}\fontshape{#5}%
  \selectfont}%
\fi\endgroup%
{\renewcommand{\dashlinestretch}{30}
\begin{picture}(1076,1089)(0,-10)
\put(538,537){\ellipse{1060}{1060}}
\path(913,912)(913,162)
\end{picture}
}
 } 
  \end{array}
. \]
We thus find that $\calA^V_n=(\text{chord diagrams})/(\text{$\FourT$
and $\FI$ relations})$, as usual in the theory of Vassiliev finite type
invariants of knots.

The Fundamental Theorem of Finite Type Invariants, due to
Kontsevich~\cite{Kontsevich:Vassiliev}, asserts that (at least over
$\bbQ$) this is indeed all: For every $W\in(\calA^V_n)^\star$ there exists
a type $n$ invariant $I$ with $W=\partial_V^nI$. In other words, every
once-integrable weight system is fully integrable.

\section{The Goussarov definition on its own} \label{sec:Goussarov} The
purpose of this section is to tell the parallel story for the Goussarov
theory of finite type invariants. Much of the mathematical content of
this section is independent of that of the previous one. But we choose
not to repeat the formal parts of the story, and to concentrate only on
the ``new stuff''.  Thus this section cannot be read independently.

\begin{figure}[ht!]
  \begin{center}
    \includegraphics[width=1.4in]{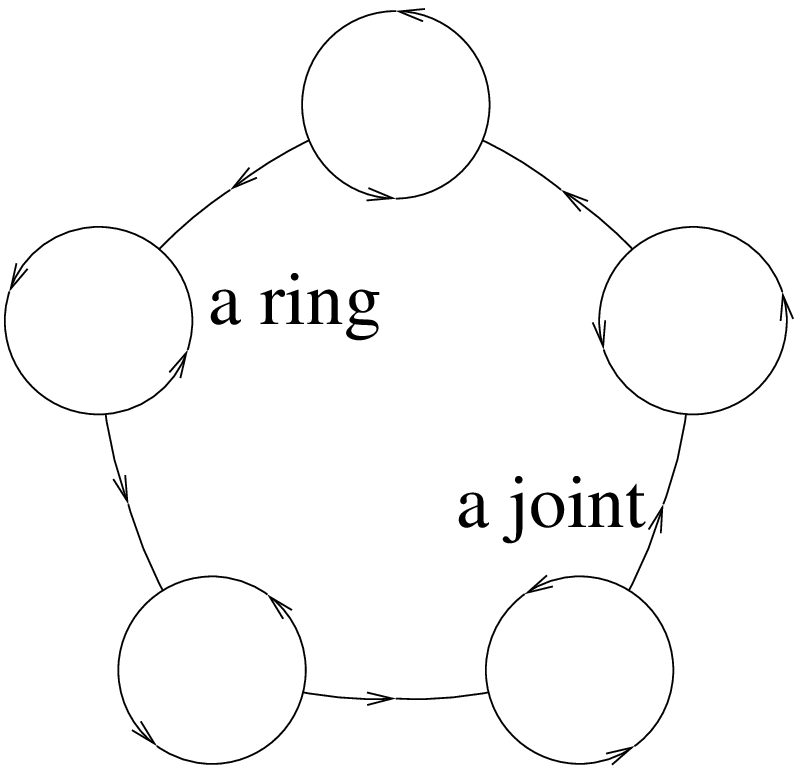}
  \end{center}
  \caption{A 5-bracelet} \label{fig:5Bracelet}
  \end{figure}

In the Goussarov theory, what replaces the space $\calK_n^V$ of formal
linear combinations of $n$-singular knots (modulo differentiability) is
the space $\calK^G_n$ of formal linear combinations of knotted
$n$-bracelets (modulo differentiability, defined later). A knotted
$n$-bracelet is an embedding up to isotopy in $\bbR^3$ of an
$n$-bracelet: a directed graph made of $n$ rings and $n$ joints. An
example of a 5-bracelet appears in
Figure~\ref{fig:5Bracelet}. Figure~\ref{fig:Rerouting} on
page~\pageref{fig:Rerouting} can be made into an example of a knotted
2-bracelet by turning the dashed lines into solid lines and adding
orientations in an appropriate manner.

The replacement for $\delta^V$ is the map
$\delta^G=\delta^G_{n+1}:\calK^G_{n+1}\to\calK^G_n$ defined by
\begin{equation} \label{eq:deltaG}
  
  \setlength{\unitlength}{0.5\standardunitlength}
  \begin{array}{c}
    {\begingroup\makeatletter\ifx\SetFigFont\undefined%
\gdef\SetFigFont#1#2#3#4#5{%
  \reset@font\fontsize{#1}{#2pt}%
  \fontfamily{#3}\fontseries{#4}\fontshape{#5}%
  \selectfont}%
\fi\endgroup%
{\renewcommand{\dashlinestretch}{30}
\begin{picture}(10662,950)(0,-10)
\put(2100.000,467.000){\arc{900.000}{1.5708}{4.7124}}
\path(1977.043,3.526)(2100.000,17.000)(1985.168,62.974)
\put(2100.000,467.000){\arc{900.000}{4.7124}{7.8540}}
\path(2222.957,930.474)(2100.000,917.000)(2214.832,871.026)
\path(2550,467)(3000,467)
\path(2880.000,437.000)(3000.000,467.000)(2880.000,497.000)
\path(750,467)(1200,467)
\path(1080.000,437.000)(1200.000,467.000)(1080.000,497.000)
\path(1200,467)(1650,467)
\path(3000,467)(3450,467)
\put(7575,392){\makebox(0,0)[lb]{\smash{{\mathmode{-}}}}}
\put(6000.000,467.000){\arc{900.000}{6.2832}{9.4248}}
\put(9300.000,467.000){\arc{900.000}{3.1416}{6.2832}}
\path(6450,467)(6900,467)
\path(6780.000,437.000)(6900.000,467.000)(6780.000,497.000)
\path(4650,467)(5100,467)
\path(4980.000,437.000)(5100.000,467.000)(4980.000,497.000)
\path(5100,467)(5550,467)
\path(6900,467)(7350,467)
\path(9750,467)(10200,467)
\path(10080.000,437.000)(10200.000,467.000)(10080.000,497.000)
\path(7950,467)(8400,467)
\path(8280.000,437.000)(8400.000,467.000)(8280.000,497.000)
\path(8400,467)(8850,467)
\path(10200,467)(10650,467)
\put(0,392){\makebox(0,0)[lb]{\smash{{\mathmode{\delta^G:}}}}}
\put(3675,392){\makebox(0,0)[lb]{\smash{{\mathmode{\mapsto}}}}}
\end{picture}
} } 
  \end{array}
.
\end{equation}
The differentiability relation is the minimal relation which makes
$\delta^G$ well defined:
\[ 
  \setlength{\unitlength}{0.5\standardunitlength}
  \begin{array}{c}
    {\begingroup\makeatletter\ifx\SetFigFont\undefined%
\gdef\SetFigFont#1#2#3#4#5{%
  \reset@font\fontsize{#1}{#2pt}%
  \fontfamily{#3}\fontseries{#4}\fontshape{#5}%
  \selectfont}%
\fi\endgroup%
{\renewcommand{\dashlinestretch}{30}
\begin{picture}(8718,1539)(0,-10)
\put(309.000,762.000){\arc{600.000}{3.1416}{6.2832}}
\path(4.218,885.601)(9.000,762.000)(62.946,873.310)
\put(309.000,762.000){\arc{600.000}{6.2832}{9.4248}}
\path(613.782,638.399)(609.000,762.000)(555.054,650.690)
\path(309,1062)(309,1512)
\path(339.000,1392.000)(309.000,1512.000)(279.000,1392.000)
\path(309,462)(309,12)
\put(1209.000,762.000){\arc{600.000}{4.7124}{7.8540}}
\path(1209,462)(1209,12)
\path(1209,1062)(1209,1512)
\path(1239.000,1392.000)(1209.000,1512.000)(1179.000,1392.000)
\put(2709.000,762.000){\arc{600.000}{3.1416}{6.2832}}
\path(2404.218,885.601)(2409.000,762.000)(2462.946,873.310)
\put(2709.000,762.000){\arc{600.000}{6.2832}{9.4248}}
\path(3013.782,638.399)(3009.000,762.000)(2955.054,650.690)
\path(2709,1062)(2709,1512)
\path(2739.000,1392.000)(2709.000,1512.000)(2679.000,1392.000)
\path(2709,462)(2709,12)
\put(3609.000,762.000){\arc{600.000}{1.5708}{4.7124}}
\path(3609,462)(3609,12)
\path(3609,1062)(3609,1512)
\path(3639.000,1392.000)(3609.000,1512.000)(3579.000,1392.000)
\put(5109.000,762.000){\arc{600.000}{4.7124}{7.8540}}
\path(5109,462)(5109,12)
\path(5109,1062)(5109,1512)
\path(5139.000,1392.000)(5109.000,1512.000)(5079.000,1392.000)
\put(6009.000,762.000){\arc{600.000}{3.1416}{6.2832}}
\path(5704.218,885.601)(5709.000,762.000)(5762.946,873.310)
\put(6009.000,762.000){\arc{600.000}{6.2832}{9.4248}}
\path(6313.782,638.399)(6309.000,762.000)(6255.054,650.690)
\path(6009,1062)(6009,1512)
\path(6039.000,1392.000)(6009.000,1512.000)(5979.000,1392.000)
\path(6009,462)(6009,12)
\put(7509.000,762.000){\arc{600.000}{1.5708}{4.7124}}
\path(7509,462)(7509,12)
\path(7509,1062)(7509,1512)
\path(7539.000,1392.000)(7509.000,1512.000)(7479.000,1392.000)
\put(8409.000,762.000){\arc{600.000}{3.1416}{6.2832}}
\path(8104.218,885.601)(8109.000,762.000)(8162.946,873.310)
\put(8409.000,762.000){\arc{600.000}{6.2832}{9.4248}}
\path(8713.782,638.399)(8709.000,762.000)(8655.054,650.690)
\path(8409,1062)(8409,1512)
\path(8439.000,1392.000)(8409.000,1512.000)(8379.000,1392.000)
\path(8409,462)(8409,12)
\put(1959,687){\makebox(0,0)[b]{\smash{{\mathmode{-}}}}}
\put(4359,687){\makebox(0,0)[b]{\smash{{\mathmode{=}}}}}
\put(6759,687){\makebox(0,0)[b]{\smash{{\mathmode{-}}}}}
\end{picture}
} } 
  \end{array}
. \]
We let the derivative $\partial_G$ be the adjoint of $\delta^G$, and
just as in the Vassiliev theory, we can now define finite type invariants:

\begin{definition} A knot invariant $I$ (equivalently, a linear
functional on $\calK=\calK^G_0$) is of Goussarov type $n$ if
its $(n+1)$-st (Goussarov style) derivative vanishes, that is, if
$(\partial_G)^{n+1}I\equiv 0$.
\end{definition}

Just like in the Vassiliev case, we have ladders
\begin{equation} \label{eq:Gladders}
    \def\neg{\hspace{-13pt}}
    \begin{array}{cccccccccr}
      \ldots\longrightarrow\neg &
        \calK^G_{n+1} &
        \neg\stackrel{\delta^G}{\longrightarrow}\neg &
        \calK^G_n &
        \neg\stackrel{\delta^G}{\longrightarrow} &
        \calK^G_{n-1} &
        \longrightarrow\ldots &
        \stackrel{\delta^G}{\longrightarrow} &
        \calK^G_0=\calK \\
      &&&&&&&& \\
      \ldots\longleftarrow\neg &
        (\calK^G_{n+1})^\star &
        \neg\stackrel{\partial_G}{\longleftarrow}\neg &
        (\calK^G_n)^\star &
        \neg\stackrel{\partial_G}{\longleftarrow} &
        (\calK^G_{n-1})^\star &
        \longleftarrow\ldots &
        \stackrel{\partial_G}{\longleftarrow} &
        (\calK^G_0)^\star=\calK^\star. \\
      & \inup\rule{0pt}{16pt}     & & \inup &
        & & & & \inup\  \\
      & \partial_G^{n+1}I\equiv 0   & & \partial_G^nI=W &
        & & & & I\
    \end{array}
  \end{equation}
For the same reasons as in the Vassiliev case we are lead to be interested
in the space $\calK^G_n/\delta^G\calK^G_{n+1}$. This is the space on
which ``(Goussarov style) weight systems'' are defined, and it is the
parallel of the space of chord diagrams in the Vassiliev case:

\begin{proposition} \label{prop:Gsymbols}
The space $\calK^G_n/\delta^G\calK^G_{n+1}$ is canonically isomorphic
to the space $\calD^G_n$ of formal linear combinations of ``cyclically
ordered $n$-component links'', which are simply $n$-component links
along with a cyclic order on their components.
\end{proposition}

\begin{proof} Dividing by $\delta^G\calK^G_{n+1}$ is the same as imposing
the equality
\[ 
  \setlength{\unitlength}{0.5\standardunitlength}
  \begin{array}{c}
    {\begingroup\makeatletter\ifx\SetFigFont\undefined%
\gdef\SetFigFont#1#2#3#4#5{%
  \reset@font\fontsize{#1}{#2pt}%
  \fontfamily{#3}\fontseries{#4}\fontshape{#5}%
  \selectfont}%
\fi\endgroup%
{\renewcommand{\dashlinestretch}{30}
\begin{picture}(6024,933)(0,-10)
\put(1362.000,459.000){\arc{900.000}{6.2832}{9.4248}}
\put(4662.000,459.000){\arc{900.000}{3.1416}{6.2832}}
\path(1812,459)(2262,459)
\path(2142.000,429.000)(2262.000,459.000)(2142.000,489.000)
\path(12,459)(462,459)
\path(342.000,429.000)(462.000,459.000)(342.000,489.000)
\path(462,459)(912,459)
\path(2262,459)(2712,459)
\path(5112,459)(5562,459)
\path(5442.000,429.000)(5562.000,459.000)(5442.000,489.000)
\path(3312,459)(3762,459)
\path(3642.000,429.000)(3762.000,459.000)(3642.000,489.000)
\path(3762,459)(4212,459)
\path(5562,459)(6012,459)
\put(2937,384){\makebox(0,0)[lb]{\smash{{\mathmode{=}}}}}
\end{picture}
} } 
  \end{array}
. \]
In English, this equality reads ``it doesn't matter how joints are
embedded, they can be moved modulo $\delta^G\calK^G_{n+1}$''. So what
remains modulo $\delta^G\calK^G_{n+1}$ is just the manner in which
the rings are knotted. But this is precisely a cyclically ordered
$n$-component link.
\end{proof}

In the case of the Vassiliev theory, we saw that
$\calK^V_n/(\ker\delta^V+\im\delta^V)$ is the the famed space $\calA^V_n$
of chord diagrams modulo $\FourT$ and $\FI$ relations, whose dual is the space
of weight systems. To see what we get in the Goussarov theory, we first
have to understand $\ker\delta^G$.

Here are three families of elements in $\ker\delta^G$:
\begin{enumerate}
\item Let $B$ be a bracelet that has an `empty ring' --- a ring that bounds
an embedded disk that does not intersect any other ring or joint. Then
$B\in\ker\delta^G$. (Indeed, if a ring is empty then its two resolutions as
in Equation~\eqref{eq:deltaG} are isotopic).
\item Let $B$ be a bracelet and let $B'$ be the bracelet obtained from $B$
by reversing the orientation of one of the rings. Then
$B+B'\in\ker\delta^G$. (No words needed).
\item \hfill$
  \setlength{\unitlength}{0.7\standardunitlength}
  \begin{array}{c}
    {\begingroup\makeatletter\ifx\SetFigFont\undefined%
\gdef\SetFigFont#1#2#3#4#5{%
  \reset@font\fontsize{#1}{#2pt}%
  \fontfamily{#3}\fontseries{#4}\fontshape{#5}%
  \selectfont}%
\fi\endgroup%
{\renewcommand{\dashlinestretch}{30}
\begin{picture}(5724,658)(0,-10)
\put(762.000,321.000){\arc{600.000}{4.7124}{7.8540}}
\path(885.597,625.880)(762.000,621.000)(873.353,567.143)
\put(762.000,321.000){\arc{600.000}{1.5708}{4.7124}}
\path(638.403,16.120)(762.000,21.000)(650.647,74.857)
\put(2899.500,358.500){\arc{530.330}{4.5705}{6.4251}}
\path(2984.208,640.110)(2862.000,621.000)(2978.824,580.352)
\put(2824.500,358.500){\arc{530.330}{2.9997}{4.8543}}
\put(4999.500,283.500){\arc{530.330}{6.1413}{7.9959}}
\put(4924.500,283.500){\arc{530.330}{1.4289}{3.2835}}
\path(4839.792,1.890)(4962.000,21.000)(4845.176,61.648)
\path(12,321)(462,321)
\path(1062,321)(1512,321)
\path(2862,321)(3612,321)
\path(4212,321)(4962,321)
\path(5082.000,351.000)(4962.000,321.000)(5082.000,291.000)
\path(4962,321)(5712,321)
\path(2112,321)(2862,321)
\path(2742.000,291.000)(2862.000,321.000)(2742.000,351.000)
\put(687,246){\makebox(0,0)[lb]{\smash{{\mathmode{\scriptstyle B}}}}}
\put(2787,396){\makebox(0,0)[lb]{\smash{{\mathmode{\scriptstyle B'}}}}}
\put(4887,96){\makebox(0,0)[lb]{\smash{{\mathmode{\scriptstyle B''}}}}}
\end{picture}
} } 
  \end{array}
$\hfill\ \newline
Let $B$, $B'$ and $B''$ be bracelets related as above. (To be
specific: All rings and joints may be knotted, including the parts drawn
above. The parts not shown must be knotted in the same way for
$B$, $B'$ and $B''$. And finally, apart from orientations any two of $B$,
$B'$ and $B''$ share a ``half-ring''.) Then $B'+B''-B\in\ker\delta^G$. (No
words needed).
\end{enumerate}

Let $\lin$ be the span of these three families within $\ker\delta^G$. The
rationale for this name is that modulo $\lin$, bracelets become
``multi-linear'' in ``the span of their rings'' (with the third type
of elements, for example, becoming ``additivity relations''). Anyway,
in $\calK^G_n/\lin$ we can use this ``linearity'' repeatedly (and also some
isotopies) to subdivide the span of rings to tiny pieces that contain
very little:
\begin{equation} \label{eq:subdivision}
  \begin{array}{c}
    \includegraphics[height=0.5in]{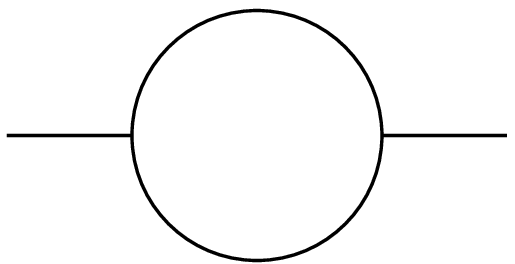}
  \end{array}
  = \ \cdots\  +
  \begin{array}{c}
    \includegraphics[height=0.5in]{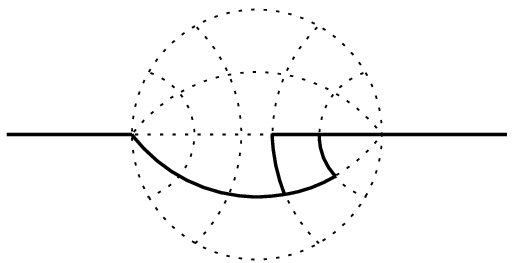}
  \end{array}
  + \ \cdots\ 
\end{equation}
Hence $\calK^G_n/\lin$ is spanned by a rather simple type of bracelets:

\noindent\parbox{3.9in}{
  \begin{definition} We say that a bracelet has simple rings if all of
    its rings bound embedded disks whose interior intersects the bracelet
    transversely and exactly once. (See an example of a simple ring on the
    right).
  \end{definition}
} \hfill $\begin{array}{c}
  \includegraphics[width=0.9in]{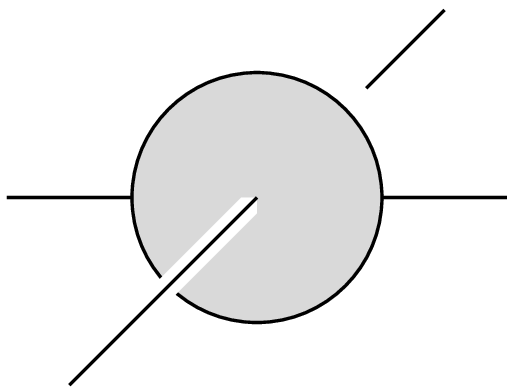}
\end{array}$

\begin{proposition} \label{prop:SimpleRings}
The space $\calK^G_n/\lin$ is spanned by bracelets with simple rings.
\end{proposition}

\begin{proof} Let $B$ be an $n$-bracelet. Find $n$ immersed disks whose
boundaries are the rings of $B$ so that there are no triple intersections
between them (this is easy; you can even arrange those $n$ disks to
have at most clasp intersections). Now subdivide all of those disks
to pieces of uniform small size as in Equation~\eqref{eq:subdivision}
(make those subdivisions sufficiently generic so that the different mesh
lines do not intersect each other and/or the joints). If the pieces are
small enough, they must be empty (and hence zero mod $\lin$) or at most
one thing may cut through any given piece.
\end{proof}

It is time for the chord diagrams of the Vassiliev theory to make their
appearance in the Goussarov theory:

\begin{proposition} For even $n$, the space $\calK^G_n/(\lin
+ \im\delta^G)$ (which is still bigger than the desired
$\calK^G_n/(\ker\delta^G + \im\delta^G)$) is isomorphic to the space of
$\frac{n}{2}$-chord diagrams. For odd $n$ the space $\calK^G_n/(\lin 
+ \im\delta^G)$ is empty.
\end{proposition}

\begin{proof}
By Proposition~\ref{prop:SimpleRings} we can reduce to bracelets
with simple rings and as in Proposition~\ref{prop:Gsymbols} we may
forget their joints. What remains is cyclically ordered $n$ component
links, each of whose components is ``simple'', meaning that it forms
a Hopf link with another component, and there's no further knotting or
linking. For odd $n$, such pairing of the components is impossible. For
even $n$ we have a cyclically ordered set of size $n$ (the components)
whose elements are paired up. This is exactly a chord diagram with $n$
vertices and $\frac{n}{2}$ chords.

\begin{floatingfigure}[l]{0.6in}
\begin{center}
  \includegraphics[width=0.6in]{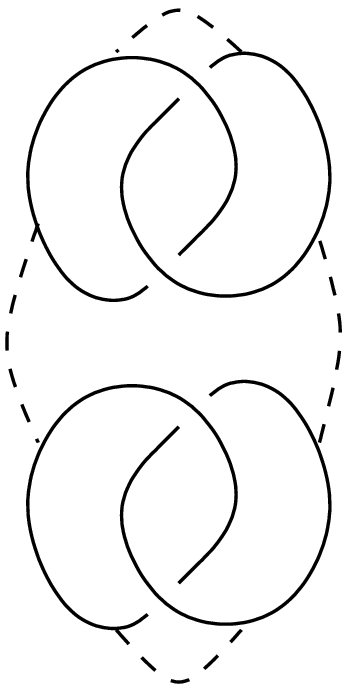}
\end{center}
\end{floatingfigure}
As an example, the figure on the left shows a bracelet with simple rings
whose corresponding chord diagram is $
  \setlength{\unitlength}{0.5\standardunitlength}
  \begin{array}{c}
    {%
\begingroup\makeatletter\ifx\SetFigFont\undefined%
\gdef\SetFigFont#1#2#3#4#5{%
  \reset@font\fontsize{#1}{#2pt}%
  \fontfamily{#3}\fontseries{#4}\fontshape{#5}%
  \selectfont}%
\fi\endgroup%
{\renewcommand{\dashlinestretch}{30}
\begin{picture}(352,363)(0,-10)
\put(176,174){\ellipse{336}{336}}
\path(26,99)(326,99)
\path(26,249)(326,249)
\end{picture}
}
 } 
  \end{array}
$. As
appropriate when moding out by $\im\delta^G$, the joints appear
``transparent''.

One still needs to show
that ``Hopf pair bracelets'' such as the one on the left, which directly
correspond to chord diagrams, do not get killed or identified with each
other by $\lin$. This can be done by noting that appropriate products
of linking numbers of rings detect Hopf pair bracelets and annihilate
$\lin$. We leave the details to the reader.
\end{proof}

There are two further families of elements in $\ker\delta^G$, the $\GFourT$
elements and the $\GFI$ elements, shown in Figure~\ref{fig:G4TFI}. We leave
it to our readers to verify that modulo $\im\delta^G$ these elements become
the $\FourT$ and the $\FI$ relations between chord diagrams:

\begin{figure}[ht!]
\[ \includegraphics[height=0.8in]{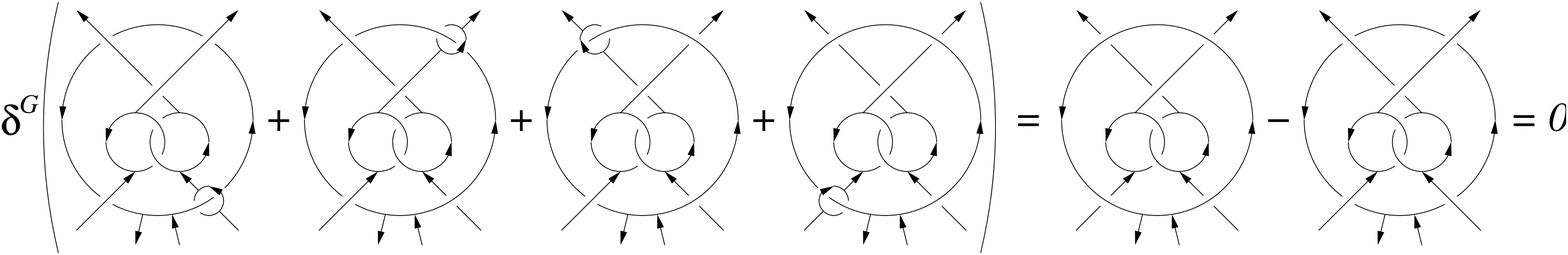} \]
\,\hfill \parbox[b]{3.9in}{
  \caption{
    The $\GFourT$ family of elements of $\ker\delta^G$ (above) and the $\GFI$
    family of elements of $\ker\delta^G$ (right).
  } \label{fig:G4TFI}
} \hfill
$\begin{array}{c}
  \includegraphics[height=0.45in]{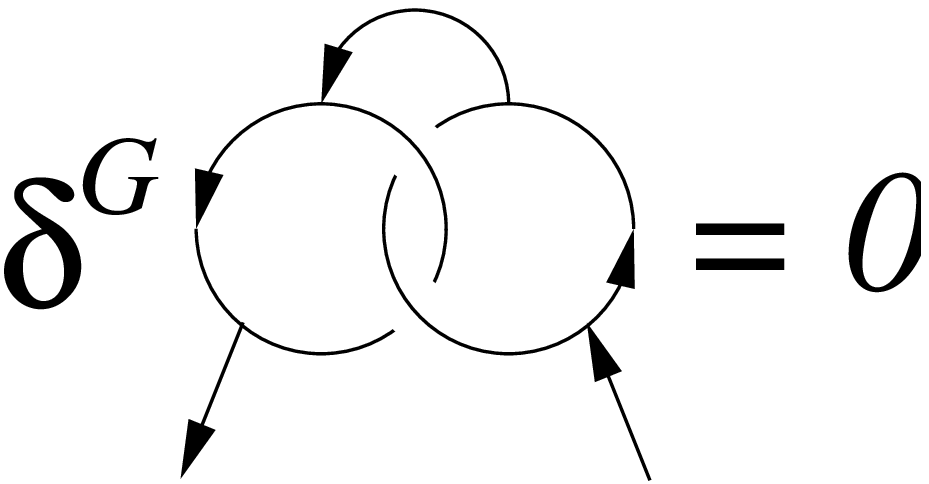}
\end{array}$ \hfill\,
\end{figure}

\begin{proposition} For even $n$, the space
$\calK^G_n/(\lin+\GFourT+\GFI+\im\delta^G)$ is isomorphic to the space
$\calA^V_{n/2}$ of the Vassiliev theory. For odd $n$ the space
$\calK^G_n/(\lin+\GFourT+\GFI+\im\delta^G)$ is empty. \qed
\end{proposition}

\begin{remark} In the light of the equivalence of the Goussarov theory
and the Vassiliev theory (shown in the next section), it is clear that
$\lin+\GFourT+\GFI+\im\delta^G=\ker\delta^G+\im\delta^G$, at
least over $\bbQ$. I do not know if $\lin+\GFourT+\GFI=\ker\delta^G$.
\end{remark}

\section{The equivalence of the two definitions} \label{sec:Equivalence}

As stated (in a slightly different form) in the introduction, the key to
the proof of Theorem~\ref{thm:main} is the (informal) equality
\[
  \delta^V\left(\begin{array}{c}
    \includegraphics[width=0.4in]{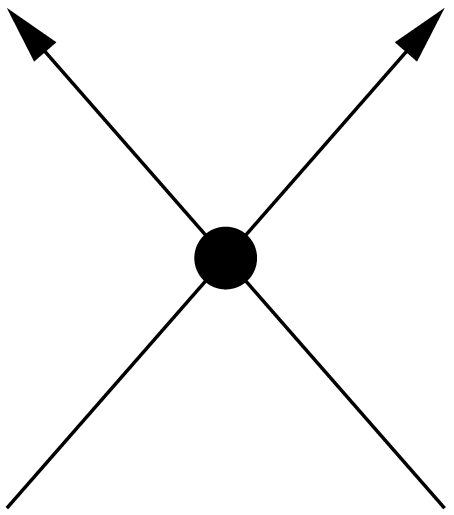}
  \end{array}\right)
  = (\delta^G)^2 \left(\begin{array}{c}
    \includegraphics[width=0.475in]{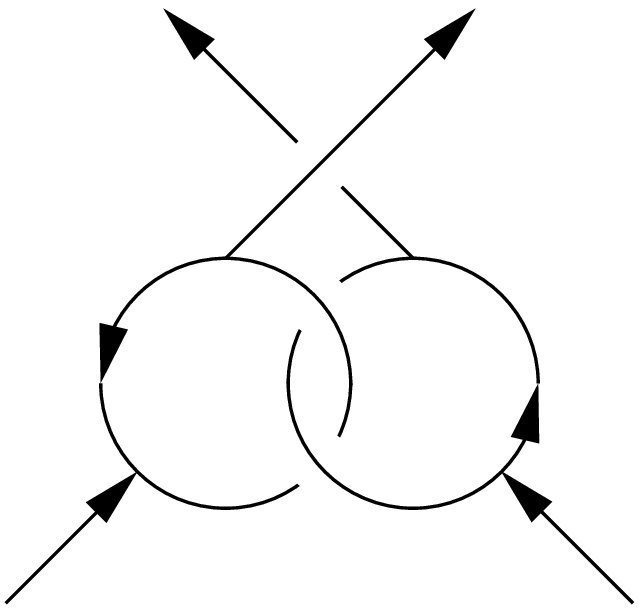}
  \end{array}\right)
\]
Let us turn this into a precise argument:

\begin{proof}[Proof of Theorem~\ref{thm:main}] An invariant $I$ is of
Vassiliev type $n$ if it vanishes on $(\delta^V)^{n+1}(\calK^V_{n+1})$
and is of Goussarov type $2n$ (respectively $2n+1$) if it
vanishes on $(\delta^G)^{2n+1}(\calK^G_{2n+1})$ (respectively
$(\delta^G)^{2n+2}(\calK^G_{2n+2})$). Thus we need to prove that
\begin{equation} \label{eq:hard}
    (\delta^G)^{2n+1}(\calK^G_{2n+1})\subset(\delta^V)^{n+1}(\calK^V_{n+1})
  \end{equation}
  and that
  \begin{equation} \label{eq:easy}
    (\delta^V)^{n+1}(\calK^V_{n+1})\subset(\delta^G)^{2n+2}(\calK^G_{2n+2}).
  \end{equation}
\parpic[r]{\includegraphics[width=1in]{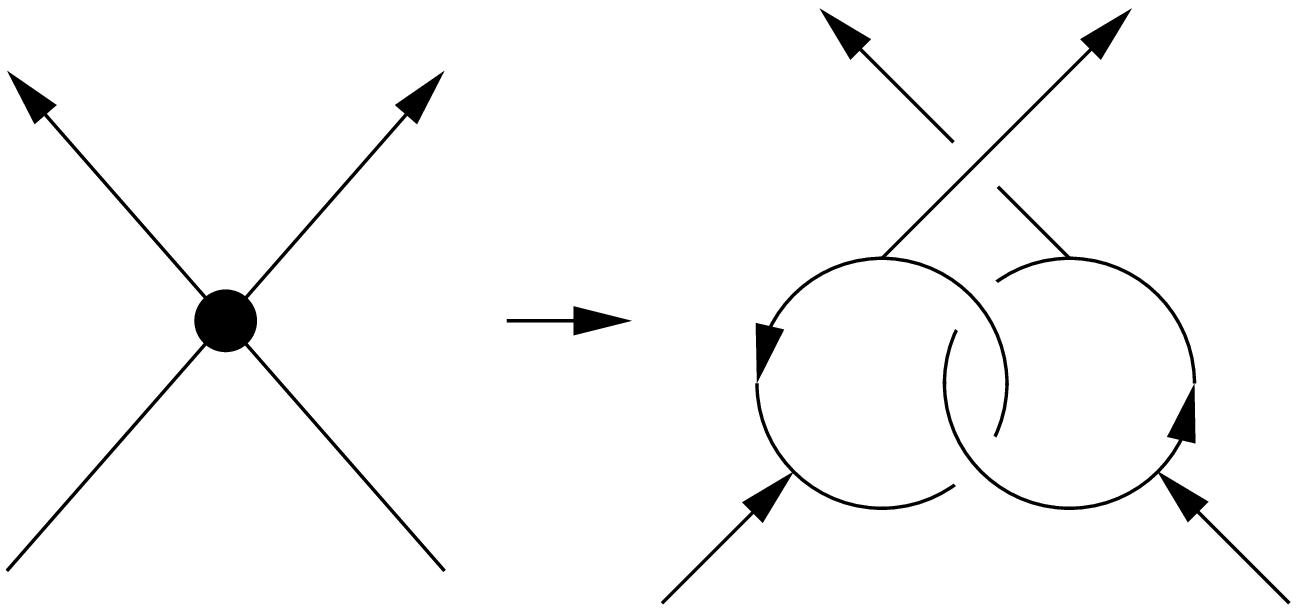}}
The easier part is the proof of~\eqref{eq:easy}. Let
$K\in\calK^V_{n+1}$ be an $(n+1)$-singular knot, and let
$B_K\in\calK^G_{2n+2}$ be the $(2n+2)$-bracelet obtained from $K$ by
replacing every singular point with a pair of rings using the rule on
the right. It is clear that
$(\delta^V)^{n+1}(K)=(\delta^G)^{2n+2}(B_K)$, and as $K$ was arbitrary,
this proves~\eqref{eq:easy}.

\parpic[l]{\includegraphics[width=0.7in]{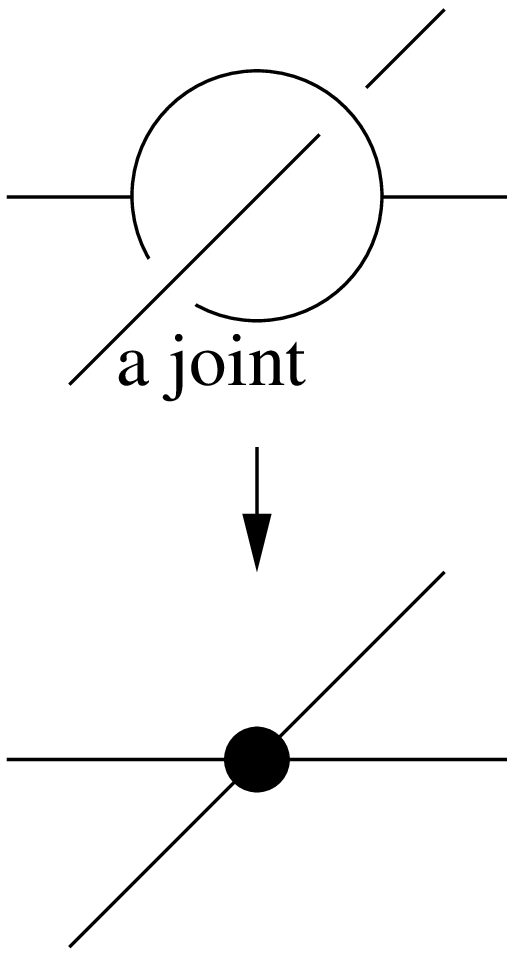}}
Let us now prove~\eqref{eq:hard}. Let $B\in\calK^G_{2n+1}$ be a
$(2n+1)$-bracelet. We need to show that $(\delta^G)^{2n+1}(B)$ is in
$(\delta^V)^{n+1}(\calK^V_{n+1})$. Clearly it does not matter if we
modify $B$ by adding to it elements in $\ker\delta^G$, so using
Proposition~\ref{prop:SimpleRings} we may assume that $B$ has simple
rings. A simple ring may loop around a joint or it may be Hopf-linked
with another simple ring. In the former case, apply the rule on the
left. In the latter case, apply the reverse of the rule in the first
half of the proof. Doing so to all rings we get a singular knot $K_B$
that has at least $n+1$ singularities (every ring in $B$ contributes
either $1$ singularity or $\frac12$ singularity, and $B$ has $2n+1$
rings). If $K_B$ has $m$ singularities (with $m\geq n+1$), we have
$(\delta^G)^{2n+1}(B) = (\delta^V)^m(K_B) \in (\delta^V)^m(\calK^V_m)
\subset (\delta^V)^{n+1}(\calK^V_{n+1})$.
\end{proof}

\Addresses


\begin{thebibliography}

\bibitem{Bar-Natan:OnVassiliev} {\bf D~Bar-Natan},
  {\em On the Vassiliev knot invariants},
  Topology {34} 423--472 (1995)

\bibitem{Conant:OnGoussarov} {\bf J~Conant},
  {\em On a theorem of Goussarov},
  Journal of Knot Theory and its Ramifications, to appear,
  {\tt arXiv:math.GT/0110057}

\bibitem{Goussarov:InterdependentModifications} {\bf M\,N~Goussarov},
  {\em Interdependent Modifications of Links and
    Invariants of Finite Degree}, 
  Topology 37 (1998) 595--602

\bibitem{Kontsevich:Vassiliev} {\bf M~Kontsevich},
  {\em Vassiliev's knot invariants},
  Adv.{} in Sov.{} Math., 16(2) (1993) 137--150

\bibitem{Stanford:FiniteType} {\bf T~Stanford},
  {\em Finite type invariants of knots, links, and graphs},
  Topology 35 (1996) 1027--1050

\end{thebibliography}
\end{document}